\documentclass[11pt]{article} 

\usepackage[a4paper]{geometry}
\usepackage{amsfonts,amsmath,amssymb}
\usepackage{epsfig,multicol} 
\usepackage[none]{hyphenat}
\usepackage[T1]{fontenc}
\usepackage{times}
\newtheorem{theorem}{Theorem}
 
\usepackage{mathtext}
\usepackage{multicol}
\usepackage{graphicx}
\usepackage{float}
\usepackage[belowskip=-10pt , aboveskip=5pt,justification=centering]{caption}
\usepackage{enumitem}
\usepackage{fancybox}
\usepackage{mdframed}
\usepackage{listings}
\usepackage{color}
\usepackage[usenames,dvipsnames,svgnames,table]{xcolor}

\usepackage{url}
\usepackage{enumitem}
\usepackage{fancybox}
\usepackage{mdframed}
\usepackage{listings}
\usepackage{fancyhdr}
\fancypagestyle{plain}
\usepackage{imakeidx}
\makeindex
\usepackage{sectsty}
\usepackage{lipsum}
\usepackage[hang]{footmisc}
\usepackage[autostyle]{csquotes}
\usepackage[margin=10pt,font=small]{caption}
\usepackage{url}
\usepackage[intoc,refpage]{nomencl}
\makenomenclature

\usepackage[pdftex]{hyperref}
\hypersetup{hyperindex=true,hypertexnames=true,
            hyperfootnotes=true,colorlinks=true,
            bookmarks=true,
            linkcolor=black,citecolor=black,
            filecolor=black,urlcolor=black}

\begin{document}

\title {\Huge{\textbf{New models for deformations:}} \\ {\LARGE{\textbf{ Linear Distortion and the failure of rank-one convexity}} }}
\author{\large{\textbf{ Seyed Mohsen Hashemi \:\:\& \: Gaven J. Martin\footnotemark[1]}}}
\footnotetext[1]{Work of both authors partially supported by the New Zealand Marsden Fund.\\Institute for Advanced Study, Massey University, Auckland, New Zealand\\
email: g.j.martin@massey.ac.nz}
\date{\textbf{April 2020}}
\maketitle

\begin{abstract}
\noindent In this article, we discuss new models for static nonlinear deformations via scale-invariant conformal energy functionals based on the linear distortion. In particular, we give examples to show that, despite equicontinuity estimates giving compactness, minimising sequences will have strictly lower energy than their limit, and that this energy gap can be quite large. We do this by showing that Iwaniec's theorem on the failure of rank-one convexity for the linear distortion of a specific family of linear mappings, is actually generic and we subsequently identify the optimal rank-one direction to deform a linear map to maximally decrease its distortion. 
\end{abstract}
\section{Introduction}

\subsection{Aspects of material science}

The use of nonlinear elasticity to describe martensitic transformations and their microstructure starts with J.M. Ball and R.D. James (1987) \cite{BJ}, following work of many authors applying nonlinear elasticity to crystals, especially J.L. Ericksen. There is a ``linearized'' version of the theory due to Khachaturyan and Roitburd. The lecture notes of Ball \cite{Ball} give a very good overview of the subject at an introductory level. In fact, one of the most studied phase transitions in nonlinear materials science is the Martensitic transition, diffusionless structural transition from a high symmetry to a lower symmetry crystallographic phase. This transition is often induced by changing temperature or stress and shows athermal character and proceeds intermitte-\\ntly as a sequence of avalanches ultimately producing a complex multiscale microstructure.\\
Recently the theory of mappings of finite distortion has emerged and shown promise in modelling various aspects of nonlinear materials science with interesting associated extremal problems, see for example \cite{AIM1,AIMO,HK2,IKO1,IM1,IMO,IO1,IO2,MY1}. Many of the distortion functionals studied in these cases are lower semicontinuous and extremals are quite regular - often diffeomorphisms. However, the images above suggest that this cannot be the case for the sort of phase transitions occuring in these structural transformations. 

\begin{figure}[ht]
         \begin{center}
         \includegraphics[scale=0.48]{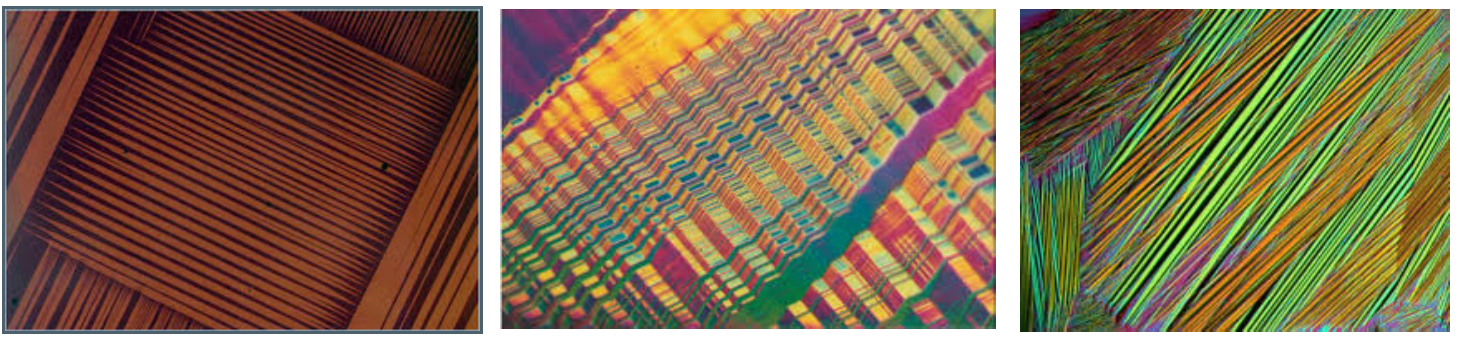}
         \caption{These images are from Chu and James's experiments in Cu-Al-Ni single crystals,\\ and M. Morin \cite{Ball}}
         \end{center}
\end{figure}

Thus we analyse a slightly different (but still well known)  distortion functional where we loose polyconvexity and hope to gain models exhibiting non-smooth extrema. Here we are only concerned with static models but are well aware that the real problem is dynamical. We hope to address this elsewhere. Distortion functionals are polyconvex in two-dimensions,  so it is unlikely we will see these complicated microscale structures with purely two dimensional models using distortion. However, an interesting feature of our models is that although fundamentally three-\\dimensional, limiting processes yield two-dimensional and even one-dimensional models that at least suggest connections.


\subsection{Distortion functions} 

Distortion functionals are scale-invariant measures of  the anisotropic nature of a deformation.  We view their $L^p$-norms as stored energy functionals,  and the $L^\infty$-norm as a ``conformal energy''.  All are measures of the deviation from a conformal mapping,  and in the problems we consider, should a conformal mapping be a candidate, it will be an absolute minimiser. For the purposes of this article, we restrict two three dimensions and consider deforming a body $\Omega\subset \mathbb{R}^3$ by a  {\color{black}\textbf{ homeomorphism}} $f:\Omega\rightarrow f(\Omega)\subset \mathbb{R}^3$. Typically the assumption that a deformation is a homeomorphism is given to us by the principle of interpenetrability of matter,  see \cite{BallIPM}. It is a remarkable feature of mappings of bounded distortion, coming from modulus of continuity estimates on mappings and their inverses, that this topological condition is retained under limits.
Let $\Omega\subset \mathbb{R}^n$,  $n\geq 3$,  be a domain and $f:\Omega \longrightarrow  \mathbb{R}^n$ be a homeomorphism. Typically one assumes some regularity of the deformations in question, namely
 $$f\in W^{1,3}_{loc}(\Omega),$$ 
 the Sobolev space of mappings whose first derivatives are in $L^3$. This condition is basically the minimal assumption one can make to ensure that the Jacobian determinant is locally integrable and that something like the change of variables formula might hold \cite{IM1}.

\subsection{Polyconvex distortion functions}

The distortion inequality assumes that the full differential of a mapping is controlled by the Jacobian; there is a constant $K<\infty$ such that
\begin{equation}
\| Df(x) \|^2 \leq  K J(x,f), \hskip20pt \mbox{almost every $x\in \Omega$}.
\end{equation}
Such mappings are called {\color{black}\textbf{quasiconformal}}. Deformations satisfying the distortion inequality enjoy many properties, including higher regularity $$f\in W^{1,3+\epsilon(K)}_{loc},$$  positive Jacobian $J(x,f)>0$ almost everywhere and the change of variables formula \cite{IM1}. Thus one can define the distortion function\\
\begin{equation}
\mathbf{K}(x,f) = \frac{\|Df(x)\|^2}{J(x,f)}
\end{equation}\\
It is not obvious, but not too hard to see that $\mathbf{K}(x,f)$ is a convex function of the minors of $Df(x)$. Thus one can consider various extremal problems, for instance, minimising the $L^p$-norms of $\mathbf{K}(x,f)$ among deformations with prescribed boundary values linking with the calculus of variations. Some novel phenomena arise in such problems including the Nitsche phenomena where \cite{AIMO} geometric obstructions preclude the existence of a minimiser. However, it is expected that minimisers are smooth when they exist (and they should exist for all $p>1$) if there is a candidate ``barrier''.

\subsection{Definition of linear distortion} 

For $x\in\Omega$, we set
\begin{equation}\label{LD}
{\mathcal{H}(x,f) = \limsup_{r\to 0} \;\frac{\max_{|h|=r}|f(x+h)-f(x)|}{\min_{|h|=r}|f(x+h)-f(x)|}}
\end{equation}
Then set 
 $\mathcal{H}(f) = {\rm ess\; sup}\{\mathcal{H}(x,f):x\in \Omega \}. $ 

If $\mathcal{H}(x,f)$ is bounded in $\Omega$,  then $f$ is $K$-quasiconformal for some $K$ bounded by a function of $\mathcal{H}$.  In particular $f$ lies in the Sobolev space $W^{1,3}_{loc}(\Omega)$ and satisfies the  {\color{black}\textbf{distortion inequality}},
$$
{ \max_{|\zeta|=1} |Df(x) \zeta| \leq \mathcal{H}(f) \min_{|\zeta|=1} |Df(x)\zeta|}.
$$
Geometrically this means that the differential $Df(x):\mathbb{R}^n \longrightarrow \mathbb{R}^n$ maps the unit sphere to an ellipsoid. The converse is also true,  but a remarkable theorem of Heinonen and Koskela (1995)  shows that one only requires the  $\mathbf{liminf}$ in (\ref{LD}) to gain quasiconformality (with the same constants).
{\color{black}\textbf{A natural question is:}} \textbf{\em{Can these scale-invariant problems and Martensite transitions be modelled by distortion functionals} ?} This relationship (if it exists at all) is predicated that at molecular scales we see piecewise linear mappings of uniformly bounded distortion as illustrated below. There are also some conjectural ideas regarding boundary values and questions such as are local scales  determined by 
{\color{black}\textbf{``interfacial energy''}} ? \\

\begin{figure}[ht]
         \centering
         \includegraphics[scale=0.6]{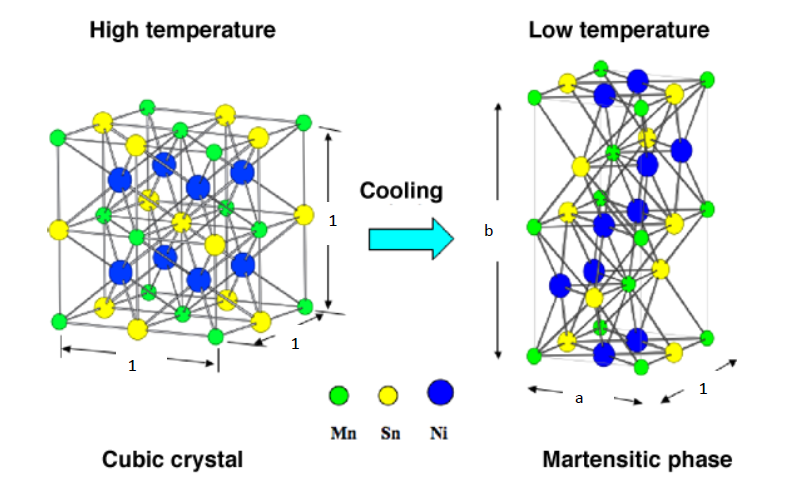}
         \caption{Martensitic phase after cooling down process. }
\end{figure}
\medskip
\begin{figure}[thp]
\centering
\begin{minipage}[b]{0.45 \textwidth}
      \centering
      \includegraphics[scale=0.42]{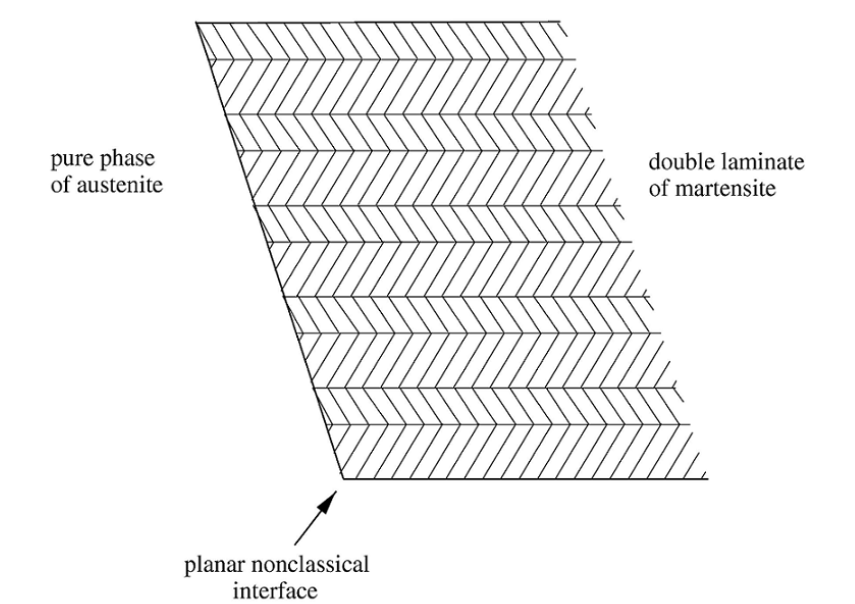}
      \caption{Nonclassical interface with\\  double laminate \cite{BJ}.}
\end{minipage}%
\hfill%
\begin{minipage}[b]{0.55 \textwidth}
      \centering
      \includegraphics[scale=0.4]{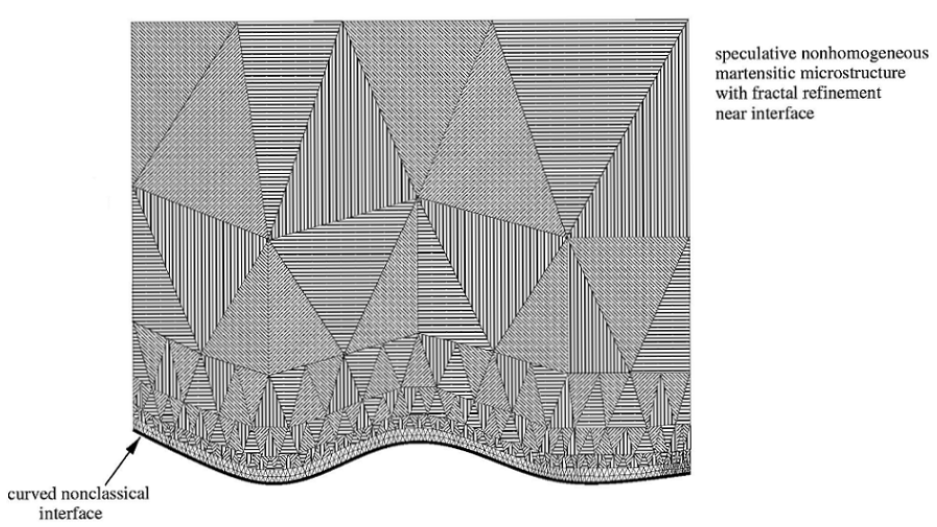}
      \caption{Nonclassical austeite-\\martensite interface \cite{BJ}.}
\end{minipage}
\end{figure}
\medskip
\section{A Model Problem}

Let $\mathbf Q=[0,1]^3$ be the unit cube in $\mathbb{R}^3$ and $A\in GL(3,\mathbb{R})$.  Consider the problem:  minimize 

\begin{equation}\label{qc} 
\int_\mathbf{ Q} \Phi[\mathcal{ K}(x,f)] \; dx \;\;\; \;\big|\;\;\;\ f:\mathbf{ Q} \longrightarrow \mathbb{R}^n \;\; \mbox{is quasiconformal.    } 
 \end{equation}

Here $\Phi$ is convex increasing and $\mathcal{K}(x,f)$ is a distortion function.   We want to discuss boundary constraints.  Thus label the faces of $\mathbf{ Q}$ as
 $$\sigma_j=[0,1]^j\times \{0,1\}\times[0,1]^{3-j}, \:\:\:\: j=1,2,3.$$
We consider imposing the following  possible constraints. Notice that quasiconvexity of the functional at (\ref{qc}) would guarantee that the minimum should be attained by the linear map $A$.
\begin{itemize}
\item $f|\partial \mathbf{ Q}=A$; in this case, typically the linear map is the minimiser, but is not unique for linear distortion unless $A$ has repeated singular values. 
\item $f(\sigma_j) = A(\sigma_j)$, in this case, we have faces being mapped to faces and typically a linear minimiser; or  
\item $ f(\sigma_j) \subset N_\epsilon(A(\sigma_j))$, $j=1,\dots, n$; in this case, we have faces being mapped within and $\epsilon$ neighbourhood of faces. Then the linear map is never a minimiser for linear distortion; or    
\item $  \partial f(\mathbf{ Q}) \subset N_\epsilon(A(\partial \mathbf{Q})$  as above. 
\end{itemize}
For the last two we take a limit as $\epsilon\to 0$. In all cases, the linear map is the minimiser for $\|\mathbf{K}(x,f)\|_p$,  $1\leq p \leq \infty$, \cite{AIMO}. In the last two cases, the linear map may be the limit of a minimising sequence but is not a minimum.
\begin{mdframed}[backgroundcolor=white]
\begin{theorem}[Linear maps never minimize]\label{def}
Let $\Omega\subset \mathbb{R}^n$ be a bounded domain and $A\in GL(3,\mathbb{R})$ with three  distinct singular values. Then there is $\delta>0$ with the following property:\\    
For every $\epsilon>0$,   there is a quasiconformal mapping $f_\epsilon:\Omega \to\mathbb{R}^n$ so that 
\begin{enumerate} 
\item $\|f_\epsilon (z)-A(z)\|_{L^{\infty}(\Omega)}<\epsilon$.
\item $$ \int_\mathbf{ Q} \Phi[\mathcal{H}(x,f)] < \int_{\bf Q} \Phi[\mathcal{H}(A)] - \delta .$$
\end{enumerate} 
Here $\mathcal{H}(x,f)$  is the  linear distortion . 
\end{theorem}
\end{mdframed}
This result is true in $\mathbb{R}^n$, $n\geq 3$ as soon as there are three distinct singular values for $A$.  This result exhibits the failure of lower semi-continuity. {\color{black}\textbf{An intriguing question is how big can the jump be?}}
\section{The Linear Distortion}
 
For many years - perhaps since the inception of the higher dimensional theory of quasiconformal mappings -  it was assumed that, like other measures of distortion for mappings such as the inner and outer distortion, the linear distortion was lower semi-continuous; that is {\color{black}\textbf{the distortion of a limit is never more  than the limit of the distortions}}. Tadeusz Iwaniec gave a striking example that refuted this belief in 1998. The key element of his construction is that the linear distortion function fails to be rank-one convex in dimension higher than $2$. 
\begin{mdframed}[backgroundcolor=white]
\begin{theorem}[Iwaniec \cite{Iw}]
 For each dimension  $n\geq 3$ and dilatation $H>1$, there exists a sequence $\{f_j\}_{j=1}^{\infty}$ of $H$-quasiconformal mappings of $\mathbb{R}^n$ converging uniformly to a quasico-\\nformal linear map
 $f:\mathbb{R}^n \longrightarrow \mathbb{R}^n$ whose dilatation is greater than $H$.
\end{theorem}
\end{mdframed}
That the dimension is greater than two is necessary of course.  Understanding the deep connectio-\\ns between problems of semi-continuity,  the various notions of convexity of functionals and materials science,  Iwaniec realised the question could be reduced to deciding the rank-one convexity of the linear distortion functional.  He gave a specific family of examples establishing this theorem.  Our results show that this failure is a generic feature.  We now discuss the proof of this.
\subsection{Rank-One Convexity }

A function  $\mathcal F:U\subset \mathbb{R}^{n\times n}\longrightarrow \mathbb{R}$, defined on an open set $U$ of $n\times n$ matrices, is said to be {\color{black}\textbf{rank-one convex}} at $A_0\in U$ if  for every rank-one matrix $X\in \mathbb{R}^{n\times n}$ the function
$$ t \mapsto \mathcal{F}(A_0+tX)$$
is convex near $t=0$. 
Recall that $X$ has rank one if and only if it is the tensor product of two vectors $X=u\otimes v$, where $u,v\in \mathbb{R}^n$.\\  
The Hadamard jump condition asserts that the piecewise linear function $F=\{A,B\}$ is continuous across the interface if and only if $A-B = u\otimes v$,  and the normal to the interface is $v$.
 \begin{figure}[ht]
         \centering
         \includegraphics[scale=1]{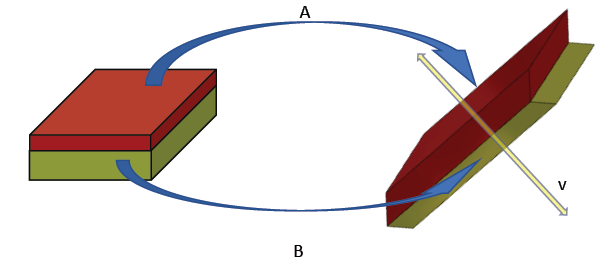}
         \caption{Rank-one connected planar interface}
\end{figure}

\bigskip

It is an exercise to show that the linear distortion functional is $\mathcal{H}:GL^{+}(3,\mathbb{R})\to \mathbb{R}$ is defined by 
$$ \mathcal{H}(A) = \frac{\max_{|\zeta|=1} |A\zeta|}{\min_{|\zeta|=1} |A\zeta| }  = \sqrt{\frac{\lambda_n}{\lambda_1}},$$
 where $\lambda_1\leq \lambda_2 \leq   \lambda_3$ are the ordered eigenvalues of $A^tA$. Iwaniec showed (by hand!) that $\mathcal{H}$ is not rank-one convex at $A_c$ for $c>1$,  where
$$ A_c= \left[\begin{array}{ccc}1&0&0 \\ 0 & c & 0 \\ 0 & 0& c^2 \end{array}\right].$$ 
Note $\mathcal{H}(A_c)=c^2$. He chose two vectors $u=(1,b,c)$ and $v=(1,-b,c)$. If $X=u\otimes v$, then he proves that 
$$H(A+t_+X)=H(A+t_-X)=\sqrt{\frac{\lambda_{max}}{\lambda_{min}}}=c^2 -2 t^2+{\rm higher\: powers\: of}\: t,$$
where $\lambda_{max}$ and $\lambda_{min}$ are the biggest and the smallest eigenvalues of matrix $(A+tX)^T(A+tX).$ It becomes obvious that 
$H(A+tX)<H(A)$ for $t$ sufficiently small. Our first result gives the best rank one direction to deform a linear mapping.

 \subsection{Optimal Rank-one Directions}

\begin{mdframed}[backgroundcolor=white]
\begin{theorem}\label{or1} Given $A\in GL^{+}(3,\mathbb{R})$ with distinct singular values.   Then  up to sign,  there is a unique rank one matrix    $B_0 = u_0\otimes v_0$,    $\|u_0\|=\|v_0\|=1$, with the following properties. 
\begin{itemize}
\item $ \frac{d}{dt}\Big|_{t=0} \;\mathcal{H}(A+tB_0) =0, \;\;\; {\rm and} \;\;\;  \frac{d^2}{dt^2}\Big|_{t=0} \;\mathcal{H}(A+tB_0) < 0,$ 
\item $ \frac{d^2}{dt^2}\Big|_{t=0} \;\mathcal{H}(A+tB_0)  \leq   \frac{d^2}{dt^2}\Big|_{t=0} \;\mathcal{H}(A+tB), $
\end{itemize}
 for all $B = u\otimes v$,    $\|u\|=\|v\|=1$,  with $\frac{d}{dt}\Big|_{t=0} \mathcal{H}(A+tB) = 0$
\end{theorem}
\end{mdframed} 
Next,  given $A$ and $B_0={\bf u}\otimes {\bf v}$ as above there is an interval $[t_-,t_+]$ on which the function $\mathcal{H}(A+tB_0)$ is smooth and convex,  and 
 \[ \mathcal{H}(A+t_- B_0)\;\; \& \;\; \mathcal{H}(A+t_+ B_0)  < \mathcal{H}(A) \]
 One can now construct a sequence converging to the linear mapping with strictly smaller linear distortion as follows.  First construct the periodic sawtooth function $a(t)$ whose graph we illustrate below.
 \begin{figure}[ht]
         \centering
         \includegraphics[scale=0.6]{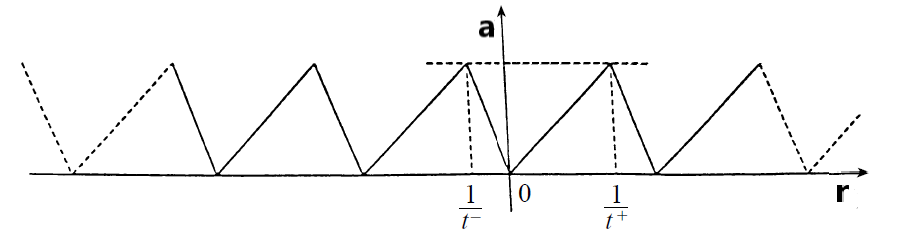}
         \caption{Sawthooth function a(t).}
\end{figure}

\bigskip

Then define function $f_j: \mathbb{R}^n \longrightarrow \mathbb{R}^n$ with equation 
$$ f_j(x) = Ax + \frac{1}{j} \; a(j {\bf u} \cdot x) {\bf v}, $$
 where $j=1,2,3,\dotsm.$ and $a$ is a piecewise linear function on $\mathbb{R}$ such that
$$ a(r)=\begin{cases}  
t^- r & \text{if \:\:$\frac{1}{t^-}\leq r \leq 0 $} \\ t^+ r & \text{if\:\: $0 \leq r \leq \frac{1}{t^+}.$}
\end{cases}  $$
 The functions $$ f_j(x) = Ax + \frac{1}{j} \; a(j {\bf u} \cdot x) {\bf v} $$  have differential which assumes only two values, which are independent of $j$:
\[ Df_j(x) = A+a'(j {\bf u}\cdot x) {\bf u}\otimes {\bf v} = A+ t_\pm B_0 \]
Then 
\[ \mathcal{H}(x,f_j) = \mathcal{H}(Df_j(x)) \leq \max \big\{ \mathcal{H}(A+ t_\pm B) \big\} < \mathcal{H}(A) \]
Effectively the piecewise linear mappings $f_j$ oscillate at finer and finer scales and converge uniformly to $f(x)=Ax$. 
 
\subsection{The maximum jump}
 
To address the size of the energy deficiency identified in Theorem (\ref{def}) we need to  find the maximum interval on which $\mathcal{H}(A+tB_0)$ is smooth and concave. Thus we need to identify the discriminant of the eigenvalue equation for $\mathcal{H}(A+tB_0)$  as it is the transverse crossing of the eigenvalues which implies $\mathcal{H}$ looses smoothness. 

\begin{figure}[ht]
         \centering
         \includegraphics[scale=0.75]{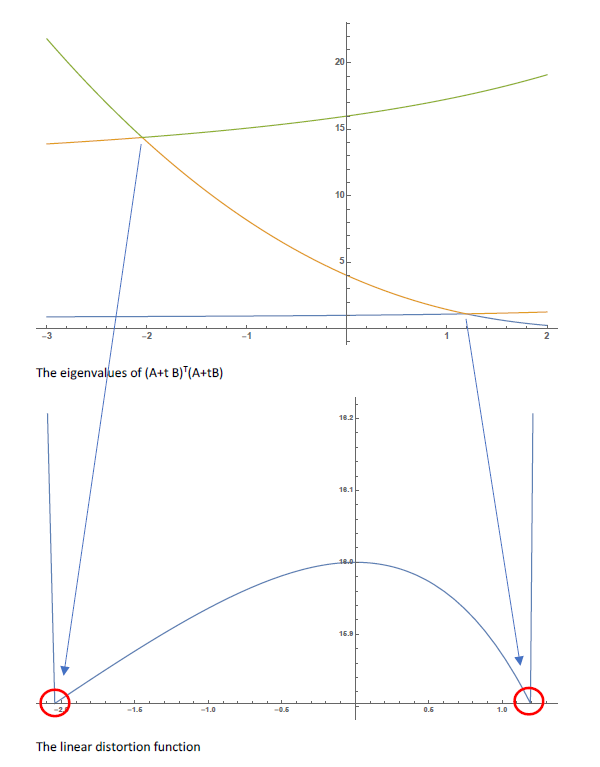}
         \caption{The function $\mathcal{H}$ is concave in the neighbourhood of $0$ that it is determined by intersection points of eigenvalues.}
\end{figure}


\section{ Optimal Directions}
 
If $A=Sing(1,\alpha,\beta)$ is a diagonal matrix, to determine the ``best'' rank one matrix $B_0=u_{0}^Tv_{0}$ with
 \begin{equation}\label{con1}
 \frac{d}{dt}\Big|_{t=0} \mathcal{H}(A+tB_0) = 0
 \end{equation}
 and that for all other rank one matrices $B$ with  $\frac{d}{dt}\Big|_{t=0} \mathcal{H}(A+tB) = 0$ we have
 \begin{equation}
 \frac{d^2}{dt^2}\Big|_{t=0} \mathcal{H}(A+t B_0) \leq  \frac{d^2}{dt^2}\Big|_{t=0} \mathcal{H}(A+t B).
 \end{equation}\vspace{0.2cm}
We need the quadratic term in the series for $ \mathcal{H}(A+tB)$ is as negative as possible.   
For $0\leq r,s \leq 1$ we put
\begin{eqnarray*}
u^T &=&  (\sqrt{1-r^2},r\cos(\theta_1),r\sin(\theta_1)), \\ v^T& = &    (\sqrt{1-s^2},s\cos(\theta_2),s\sin(\theta_2))
\end{eqnarray*}
Then $B=u^T.v$. 
Set $X=(A+tB)^T(A+tB)$ and calculate the Taylor series of the eigenvalues of $X$ to second order.  
\begin{eqnarray*}
\lambda_{min} & = &1+ \mu_1 t + \mu_2  t^2+O(t^3) \\
\lambda_{max} &=&\beta^2 +\nu_1 t + \nu_2 t^2 +O(t^3)
 \end{eqnarray*} 
Up to $O\left(t^3\right)$, 
\begin{eqnarray*}  \mathcal{H}(A+tB) & = & \beta +t \frac{\left(\eta_1-\beta ^2 \mu_1\right)}{2 \beta } \\ 
&& -\frac{t^2 \left(\beta ^4 \left(4 \mu_2-3 \mu_1^2\right)+2 \beta ^2 \eta_1 \mu_1-4 \beta ^2 \eta_2+\eta_1^2\right)}{8 \beta ^3}
 \end{eqnarray*} 
However $\mu_i,\nu_i$ are algebraic/trigonometric polynomials in four variables (six if you count $\alpha$ and $\beta$).  Then
\begin{eqnarray*}  \mathcal{H}(A+tB) & = & \beta +t \frac{\left(\eta_1-\beta ^2 \mu_1\right)}{2 \beta } \\ 
&& -\frac{t^2 \left(\beta ^4 \left(4 \mu_2-3 \mu_1^2\right)+2 \beta ^2 \eta_1 \mu_1-4 \beta ^2 \eta_2+\eta_1^2\right)}{8 \beta ^3}
 \end{eqnarray*} 
and (\ref{con1}) gives us the constraint 
\[ 0\;=\;\eta_1-\beta ^2 \mu_1 
\]
 and subject to this constraint we want to maximise
\[ {\bf Q} \; = \; \eta_2-\beta ^2 \mu_2\] 
\newpage
Now it is a long story. We use Lagrange multipliers to eliminate the variables $r$ and $s$ leading to a minimisation problem for the trigonometric algebraic function ${\bf Q}$ in two variables $\theta_1$ and $\theta_2$.  ${\bf Q}$  has various (non-obvious) symmetries which are observed by graphing a few examples and once established we see the critical points giving a minimum lie on the line $\theta_2=\pi-\theta_1$.\\

\begin{figure}[ht]
\begin{multicols}{2}
      \includegraphics[scale=0.4]{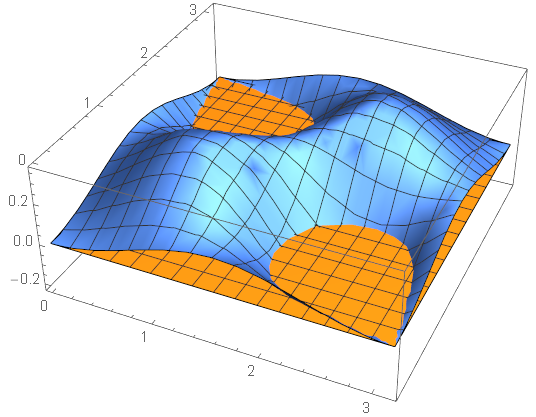}\par
       \caption{The function $\mathbf{Q}(\theta_1 ,\theta_2)$ when $A=Sing (1,2,4)$ on $[0,\pi]\times [0,\pi]$.}
      \includegraphics[scale=0.34]{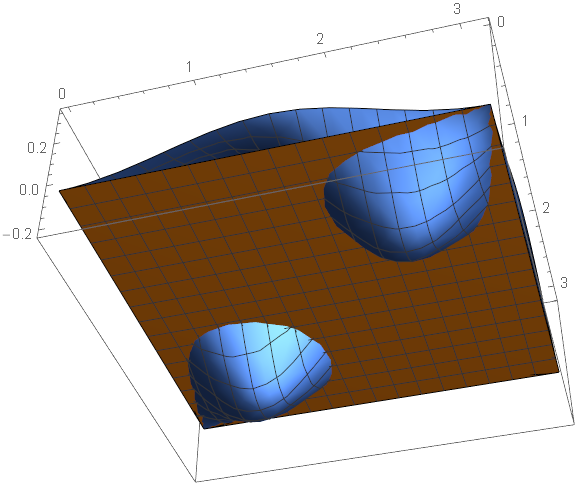}\par
      \caption{The minimums of function $\mathbf{Q}$ lie on the line $\theta_2=\pi-\theta_1$.}
\end{multicols}
\end{figure}
Now,  given the form of $u$ and $v$ we can compute the solution:  $u,v=$ 

\begin{eqnarray*}
\frac{1}{\sqrt{2} \sqrt{\alpha ^2+\alpha +\beta ^2+(\alpha -1) \beta +1}} \left(
\begin{array}{c}
(\beta -1)/\sqrt{\beta +1}   \\ \\
\pm  \sqrt{2 \alpha ^2+2 (\beta +1) \alpha +\beta ^2+1} \\ \\
(\beta -1) \sqrt{\beta }/\sqrt{\beta +1} \\
\end{array}
\right),  \end{eqnarray*} 
 
Now compute the first two terms of the series for $$\sqrt{\frac{\lambda_{max}(t)}{\lambda_{min}(t)}}.$$ To second order,
 \begin{equation*}
\mathcal{H}(A+ {\bf t}B)  =\beta - \frac{(\beta -1)^3 \beta\;  {\bf t^2} }{4 (\alpha +1) (\beta +1) (\alpha +\beta ) \left(\alpha ^2+(\alpha -1) \beta +\alpha +\beta ^2+1\right)}
 \end{equation*} 
Now it gets difficult. With this setup $\lambda_{max}$ and $\lambda_{min}$ are increasing, and $\lambda_{mid}$ is decreasing.  We have to find the largest $t_-<0$ and $t_+>0$ so the eigenvalues don't cross. Must examine the eigenvalue equation and its discriminant.
 $$ \det [(A+t B)^T(A+tB)-\lambda^2 {\bf Id} ]. $$   
One might think that
$$ \det [(A+t B)^T(A+tB)-\lambda^2 {\bf Id} ] $$
 is cubic in $\lambda^2$ and of degree $6$ in $t$.  Remarkably it is quadratic in $t$, a consequence of Jacobi's theorem. We compute the two roots in $t$ as 
 $$ \mathbf{t}= \frac{2 (\lambda -1) (\alpha +\lambda ) (\beta -\lambda )}{r^2 (\alpha  (2 \beta -\lambda -1)+(\beta -2) \lambda +\beta )-(\lambda -1) (\alpha +\beta ) \left(2 \beta  \left(r^2-1\right)+1\right)-2 (\alpha +\lambda ) (\beta -\lambda )} $$
where $\lambda$ is either square root of $\lambda^2$.  The denominator is linear in $\lambda$ and this gives us a cubic in $\lambda$ from which we can find closed form solutions for the eigenvalues.  From this we get $\lambda_{max}$ and $\lambda_{min}$ and replacing $\lambda\leftrightarrow-\lambda$ gives $\lambda_{med}$.

\section{The Eigenvalues}

If $A=Sing(1,2,10)$, then the eigenvalues of matrix  $(A+tB)^T(A+tB)$ is shown below.
\begin{figure}[ht]
         \centering
         \includegraphics[scale=0.7]{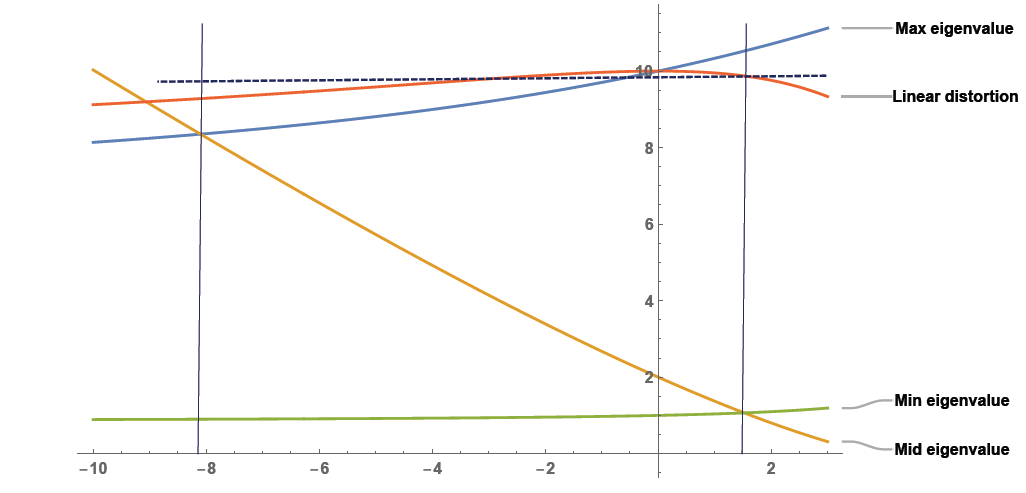}
         \caption{The eigenvalues of matrix $(A+tB)^T(A+tB)$.}
\end{figure}

\medskip

The regular branches of the  maximum and minimum eigenvalues are increasing and do not cross, while that of the middle eigenvalue is decreasing. They cross just once at asymetric points with different values of the linear distortion. The limiting jump is determined by the largest of these two values. The three eigenvalues are given by choices of $\zeta^3=-1$.

\begin{footnotesize}
\begin{eqnarray*}  -\frac{C}{3 D}-\frac{\zeta\;  \sqrt[3]{\sqrt{\left(27 A D^2-9 B C D+2 C^3\right)^2-4 \left(C^2-3 B D\right)^3}-27 A D^2+9 B C D-2 C^3}}{12 \sqrt[3]{2} D} \\ \vspace{0.3cm} +\frac{\bar\zeta \;  \left(3 B D-C^2\right)}{6\ 2^{2/3} D \sqrt[3]{\sqrt{\left(27 A D^2-9 B C D+2 C^3\right)^2-4 \left(C^2-3 B D\right)^3}-27 A D^2+9 B C D-2 C^3}} 
   \end{eqnarray*}
\end{footnotesize}
 
\begin{footnotesize}
\begin{eqnarray*} 
A&=&2 \alpha  \beta  (\beta +1) \left(\alpha ^2+(\alpha -1) \beta +\alpha +\beta ^2+1\right)-\beta  t \left(2 \alpha ^2+2 \alpha  (\alpha +4) \beta +\beta ^3+\beta ^2+\beta +1\right)  \\
B&=&t \left(2 \alpha ^2 (\beta +1)^2+\alpha  (\beta  (\beta +6)+1) (\beta +1)+\beta  (\beta  (\beta  (\beta +4)-2)+4)+1\right) \\&& -2 (\beta +1) \left(\alpha ^2+(\alpha -1) \beta +\alpha +\beta ^2+1\right) (\alpha  (\beta +1)-\beta ) \\
C&=&-2 (\beta +1) (-\alpha +\beta +1) \left(\alpha ^2+(\alpha -1) \beta +\alpha +\beta ^2+1\right)-2 (\beta +1) t \left(\alpha ^2+(\alpha -1) \beta +\alpha +\beta ^2+1\right) \\
D&=&2 (\beta +1) \left(\alpha ^2+(\alpha -1) \beta +\alpha +\beta ^2+1\right) 
\end{eqnarray*}
\end{footnotesize} 

It seems quite remarkable that the discriminant factors to give us the two crossing points,  but it is a consequence of the fact there are only two eigenvalue crossings,  yielding a quadratic factor.
\begin{footnotesize}
\begin{eqnarray*}
0&=& 2\Big(\alpha ^4 \beta ^2+2 \alpha ^4 \beta +\alpha ^4-2 \alpha ^2 \beta ^3-4 \alpha ^2 \beta ^2-2 \alpha ^2 \beta -\alpha  \beta ^5-\alpha  \beta ^4+2 \alpha  \beta ^3+2 \alpha  \beta ^2-\alpha  \beta  \\ && -\alpha +\beta ^5+\beta ^4+\beta ^2+\beta\Big)\\
&& - {\mathbf{t}_{_{\pm}}}  \Big(4 \alpha ^3 \beta ^2+8 \alpha ^3 \beta +4 \alpha ^3+\alpha ^2 \beta ^3+7 \alpha ^2 \beta ^2+7 \alpha ^2 \beta +\alpha ^2+2 \alpha  \beta ^4-4 \alpha  \beta ^3-12 \alpha  \beta ^2\\&&-4 \alpha  \beta +2 \alpha -\beta ^5-6 \beta ^4-\beta ^3-\beta ^2-6 \beta -1\Big )\\
&& -{\mathbf{t}^{2}_{_{\pm}}}\Big(-2 \alpha ^2 \beta ^2-4 \alpha ^2 \beta -2 \alpha ^2-\alpha  \beta ^3-7 \alpha  \beta ^2-7 \alpha  \beta -\alpha -\beta ^4-4 \beta ^3+2 \beta ^2-4 \beta -1\Big)
\end{eqnarray*} 
\end{footnotesize}

This gives closed form solutions for the two values where the eigenvalues cross $\mathbf{t}_{_{\pm}}$. We also have for $A=Sing(1,\alpha,\beta)$
\begin{eqnarray*}
\lim_{\beta\to\infty} \mathbf{t}_+ &=& 2(\alpha -1) \\
\lim_{\beta\to\infty} \mathbf{t}_- &=& -\infty   
\end{eqnarray*}

\section{One and Two Dimensions}

We can consider the problems in one or two dimensions by letting $\alpha,\beta \to \infty$.  The angle $\theta$ between the normal to the lamination with the principal direction $(0,0,1)$ is
\[ \cos(\theta) = \frac{(\beta -1) \sqrt{\beta }}{\sqrt{2} \sqrt{(\beta +1) \left(\alpha ^2+(\alpha -1) \beta +\alpha +\beta ^2+1\right)}} \] 

\medskip
\textbf{Quasi 1-Dimensional Case: {\em or strongly anisotropic deformation}} : 
$$\alpha/\beta\to 0\:\; \rm{we\: get\: the\: angle} \:\;\frac{\pi}{4}.$$ 
\medskip
\textbf{Quasi 2-Dimensional Case: {\em or weakly anisotropic deformation}} :\\
 If $\beta=k\alpha$, $0\leq k<1$, then
$$ \cos(\theta)=\frac{1}{\sqrt{2} \sqrt{k^2+k+1}}. $$\vspace{0.2cm}
So that, if $\beta-\alpha$ is bounded  the limiting angle is $\cos^{-1}(\frac{1}{\sqrt{6}})$.
\newpage

\begin{figure}[ht]
\begin{multicols}{2}
      \includegraphics[scale=0.44]{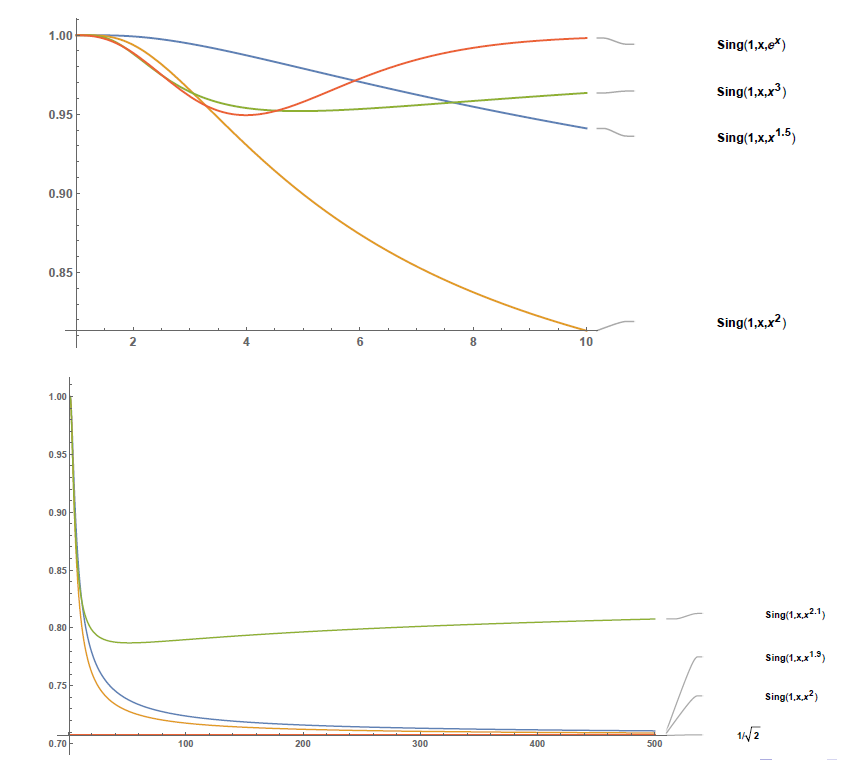}\par
       \caption{The jump - strongly anisotropic}
      \includegraphics[scale=0.35]{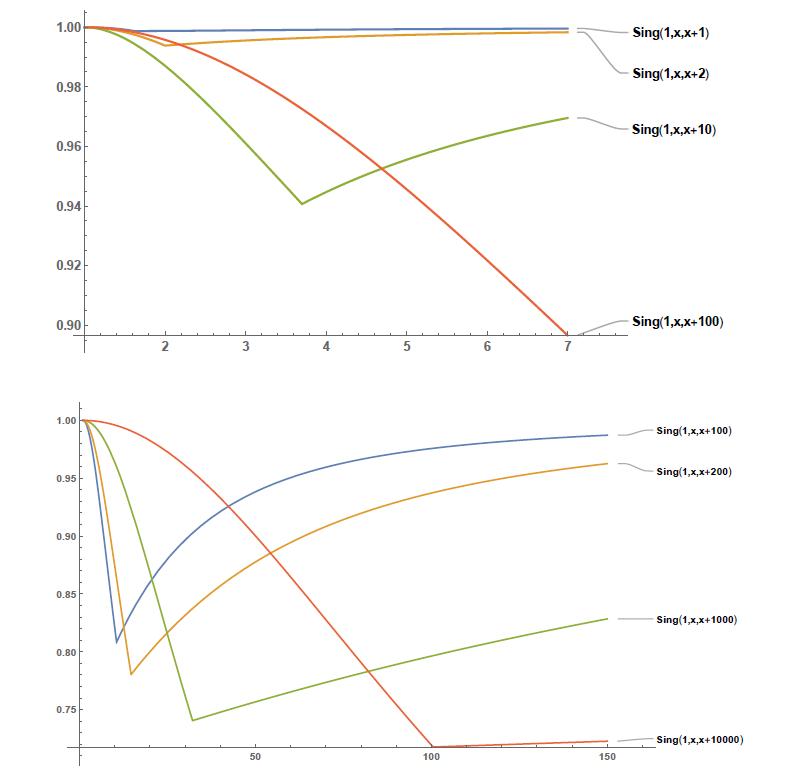}\par
      \caption{The jump - weakly anisotropic}
\end{multicols}
\end{figure}

\begin{figure}[ht]
         \centering
         \includegraphics[scale=0.8]{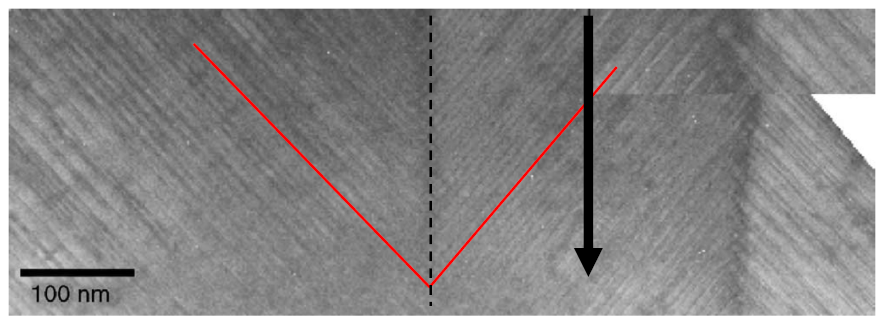}
         \caption{Austenite in strong magnetic field}
\end{figure}

\vspace{0.2cm}
 Finally,  these examples show that the energy deficiency can be as much as a factor of $\sqrt{2}$.

\begin{mdframed}[backgroundcolor=white]
 \begin{theorem}[Largest Jump] Let $c<\sqrt{2}$.  Then there is a sequence of $H$-quasiconformal mappings converging to an $H_\infty$-quasiconformal mapping (and no smaller) with 
 \[ H_\infty > c \; H\]
 \end{theorem}
\end{mdframed}

\end{document}